\documentclass{amsart}[12pt]
\usepackage{amsmath,amsfonts,amssymb,latexsym,amsthm}
\usepackage{graphicx,epsf}

\usepackage{enotez}

\usepackage{url}
\usepackage{hyperref}

\newtheorem{thm}{Theorem}[section]

\newtheorem{conv}[thm]{Convention}

\def\emp{\nothing}

\def\cc{\mathbb C}

\def\SS{\mathbb S}

\def\fq{{\mathbb F}_q}

\def\sm{\smallsetminus}

\def\la{\lambda}
\def\ga{\gamma}

\def\be{\beta}

\def\cB{\mathcal B}

\def\ssu{\subset}

\def\<{\langle}
\def\>{\rangle}

\def\rT{{\text {\rm T} } }

\newcommand{\SYT}{\operatorname{SYT}}

\def\0{{\mathbf 0}}

\def\nothing{\varnothing}

\def\.{\hskip.06cm}
\def\ts{\hskip.03cm}

\def\Ind{{\text {\rm Ind}}}

\def\bba{{\text{\bf a}}}
\def\bbb{{\text{\bf b}}}

\def\SP{{\textup{\textsf{\#P}}}}

\def\GapP{{\textup{\textsf{GapP}}}}

\def\nin{\noindent}

\begin{document}

\title{Combinatorial inequalities}

\author[Igor~Pak]{ \ Igor~Pak$^\star$}


\thanks{\thinspace ${\hspace{-.45ex}}^\star$Department of Mathematics,
UCLA, Los Angeles, CA 90095, USA.
\hskip.06cm
Email:
\hskip.06cm
\texttt{pak@math.ucla.edu}}

\vskip.2cm

\begin{abstract}
This is an expanded version of the \emph{AMS Notices} column with the
same title. The text is unchanged, but we added acknowledgements and a
large number of endnotes which provide the context and the references.
\end{abstract}

\maketitle

\vskip.25cm


{\it Combinatorics} has always been a battleground of tools and ideas.
That's why it's so hard to do, or even define.\endnote{See the discussion
in Ch.~6 of G.-C.~Rota, {\it Discrete thoughts}, Birkh\"auser, Boston, MA, 1992.
See also our blog post ``{\it What is Combinatorics?}'' at
\url{https://wp.me/p211iQ-bQ}.}  The inequalities are a
particularly interesting case study as they seem to be both the most challenging
and the least explored in Enumerative and Algebraic Combinatorics.  Here are
a few of my favorites, with some backstories.\endnote{The selection here is
somewhat arbitrary, but we tried to include some of the greatest hits which
fit a coherent narrative.  Unfortunately, many interesting results are omitted.
For example, we omitted a direct injective proof of the ad hoc inequality
\ts $E_n \cdot F_n \ge n!$, where $E_n$ is the \emph{Euler number} (the number of
alternating permutations in~$S_n$), and $F_n$ is the \emph{Fibonacci number}.  A more general
inequality \ts $e(P)\. e(P^\ast) \ge n!$ \ts for the number of linear extensions
of a $2$-dim poset and its plane dual remains open.  See A.~Morales, I.~Pak, G.~Panova, Why is $\pi < 2\phi$?,
\emph{Amer.\ Math.\ Monthly}~\textbf{125} (2018), 715--723, and our
blog post ``\emph{Fibonacci times Euler}'' at \ts \url{https://wp.me/p211iQ-lJ}.  }

\medskip

We start with {\it unimodality}\endnote{For the
general background on unimodal and log-concave sequences,
see P.~Br\"{a}nd\'{e}n, Unimodality, log-concavity, real-rootedness and beyond, in
\emph{Handbook of Enumerative Combinatorics}, CRC Press, Boca Raton, FL, 2015, 437--483;
\ts {\tt arXiv:1410.6601}.}
of binomial coefficients:

\begin{equation}\label{eq:bin}
\binom{n}{k-1} \, \le \, \binom{n}{k}\., \ \  \text{for all} \ \. 1\le k \le n/2\ts.
\end{equation}
This is both elementary and well known -- the proof is an easy calculation.
But ask yourself the following natural question: does the difference \.
$B(n,k) := \binom{n}{k} - \binom{n}{k-1}$ \. count anything interesting?  It should,
of course, right?  Imagine there is a natural injection
$$\psi: \binom{[n]}{k-1} \to  \binom{[n]}{k}
$$
from $(k-1)$-subsets to $k$-subsets of $[n]$, where $[n]:=\{1,\ldots,n\}$.
Then $B(n,k)$ can be described as the number of $k$-subsets of $[n]$
that are not in the image of~$\psi$, as good answer as any.
But how do you construct the injection~$\psi$?\endnote{There are
several proofs of this results, via injection and otherwise,
some mentioned in the main body of the paper and in the notes below.
A nice introductory overview is given in D.~Zeilberger, {\tt arXiv:1003.1273}.}


Let us sketch the construction based on the classical
\emph{reflection principle} for the \emph{ballot problem},
which goes back to the works of Bertrand and Andr\'e in~1887.
Start with a $(k-1)$-subset $X$ of~$[n]$,
and let $\ell$ be the smallest integer s.t.
$\bigl|X\cap[2\ell+1]\bigr| = \ell$. Such $\ell$ exists since $k\le n/2$.
Define
$$\psi(X)\. := \. \bigl(X \sm [2\ell+1]\bigr) \.  \cup \. \bigl([2\ell+1] \sm  X\bigr).
$$
Observe that $|\psi(X)|=k$ and check that $\psi$ is the desired injection.
This gives an answer to the original question:  $B(n,k)$
is the number of $k$-subsets $Y\ssu [n]$, s.t.
$\bigl|Y\cap [m]\bigr| \le m/2$ for all~$m$.\endnote{For details of this
proof and various generalizations see B.~E.~Sagan,
Unimodality and the Reflection Principle,
\emph{Ars Combin.}~\textbf{48} (1998), 65--72.  For the history of
the ballot problem, the reflection principle, and for further references,
see I.~Pak, History of Catalan Numbers, $\S$7, in R.~P.~Stanley,
\emph{Catalan Numbers}, Cambridge Univ.\ Press, 2015, App.~B. }

At this point you might be in disbelief in me dwelling on the easy
inequality~\eqref{eq:bin}.  Well, it only gets harder from here.
Consider, e.g., the following question: Does there
exist an injection~$\psi$ as above, s.t. $X\ssu \psi(X)$ for all \ts
$X\in \binom{[n]}{k-1}$?  We leave it to the reader as a
challenge.\endnote{One can prove this from the \emph{Hall marriage theorem}.
In fact, a stronger result holds, that the Boolean algebra
$\cB_n$ can be partitioned into symmetric saturated chains.  Use induction.
Suppose $\cB_{n-1}=\sqcup C_i$ and let $C_i':=C_i\cup \{n\}$.  Now move the
smallest element of $C_i$ to $C_i'$ and observe that $\cB_{n}=\sqcup C_i \sqcup C_i'$
as desired.  In fact, this injection can be computed in polynomial time via the elegant
``parenthesization'' algorithm.  For this construction and further references 
see C.~Greene, D.~J.~Kleitman, Strong versions of Sperner's theorem, 
\emph{J.~Combin.\ Theory Ser.~A}~\textbf{20} (1976), 80--88, and $\S$3
in C.~Greene, D.~J.~Kleitman,
Proof techniques in the theory of finite sets,
in {\it Studies in combinatorics}, MAA, 1978, 22--79.  }\label{en:GK}

There is also a curious connection to Algebraic Combinatorics: $B(n,k)=f^{(n-k,k)}$,
the dimension of the irreducible $S_n$-module corresponding to the partition
$(n-k,k)$.   To understand how this could happen, think of
both sides of~\eqref{eq:bin} as dimensions of permutation representations
of~$S_n$.  Turn both sides into vector spaces and modify $\psi$ accordingly,
to make it an $S_n$-invariant linear map.  This would make it more natural
and uniquely determined.\endnote{This is equivalent to saying that
$$
M^{(n-k,k)} \. := \. \Ind_{S_k\times S_{n-k}}^{S_n} 1 \. = \.
\SS^{(n-k,k)} \ts \oplus \. \SS^{(n-k+1,k-1)} \ts\oplus \ts \ldots \ts \oplus \. \SS^{(n-1,1)} \ts \oplus \. \SS^{(n)},
$$
where $\SS^\la$ is the irreducible $S_n$-module corresponding to~$\la$.
The uniqueness of the $S_n$-invariant injective map \ts
$\psi: M^{(n-k+1,k-1)} \to M^{(n-k,k)}$ \ts
follows from the unique decomposition in the multiplicity free case. }
As a consequence, we obtain a combinatorial
interpretation \ts $B(n,k)=\bigl|\text{SYT}(n-k,k)\bigr|$, the number of
standard Young tableaux of shape $(n-k,k)$, a happy outcome in every way.

\medskip

\nin
Consider now {\it unimodality of Gaussian coefficients:}
\begin{equation}\label{eq:qbin}
p(n,k,\ell-1) \, \le \, p(n,k,\ell)\., \ \  \text{for all} \ \. 1\le \ell \le k(n-k)/2\ts, \ \. \text{where}
\end{equation}
$p(n,k,\ell)$ is the number of integer partitions $\la\vdash\ell$ that fit
into a $k\times (n-k)$ rectangle, i.e.\ $\la$ has parts of size at most~$(n-k)$,
and has at most $k$ parts. To understand the context of this inequality, recall:
$$
\sum_{\ell=0}^{k(n-k)} \, p(n,k,\ell)\. q^\ell \, = \, \binom{n}{k}_q
\,:=\, \frac{(n!)_q}{(k!)_q \cdot \bigl((n-k)!\bigr)_q}\,, \quad \text{where} \quad (n!)_q \. := \. \prod_{i=1}^{n} \.\frac{q^i-1}{q-1}\..
$$
To connect this to~\eqref{eq:bin}, note that $\binom{n}{k}_1=\binom{n}{k}$, 
and that $\binom{n}{k}_q$ is the number of $k$-subspaces
of $\fq^n$.  In~\eqref{eq:qbin}, we view $\binom{n}{k}_q$ as a polynomial 
in~$q$ and compare its coefficients.
Now, the Schubert cell decomposition of the Grassmannian over~$\fq$, or a simple induction can be used to
give the partition interpretation.\endnote{For a full explanation of this argument see $\S$1.10 in
R.~P.~Stanley, {\it Enumerative Combinatorics}, Vol.~1 (Second ed.), Cambridge Univ.\ Press, 2012.
See also I.~Pak, When and how $n$ choose $k$, in {\it AMS DIMACS Ser.}, vol.~43 (1998), 191--238.}

The inequality~\eqref{eq:qbin} is no longer easy to prove.  Conjectured by Cayley in 1856,
it was established by Sylvester in~1878; the original paper is worth reading even if just to
see how pleased Sylvester was with his proof.  In modern language, Sylvester defined the
$\mathfrak{sl}_2(\cc)$ action on certain homogeneous polynomials and the result follows
from the highest weight theory (in its simplest form for
$\mathfrak{sl}_2$).\endnote{For a friendly exposition of this approach, 
see R.~A.~Proctor, Solution of two difficult combinatorial problems with 
linear algebra, {\it Amer.\ Math.\ Monthly}~\textbf{89} (1982), 721--734. 
A different elegant proof is given in $\S$23 of J.~Matou\v{s}ek, 
\emph{Thirty-three miniatures}, AMS, Providence, RI, 2010.  
For a more general setting, see R.~P.~Stanley, 
Unimodal sequences arising from Lie algebras,
in {\it Lecture Notes Pure Appl.\ Math.}~\textbf{57},
Dekker, New York, 1980, 127--136.}

Let's continue with the questions as we did above.  Consider the difference
$C(n,k,\ell) := p(n,k,\ell)-p(n,k,\ell-1)$.  Does $C(n,k,\ell)$ count anything interesting?
Following the pattern above, wouldn't it be natural to define some kind of nice
injection from partitions of size $(\ell-1)$ to partition of size~$\ell$,
by simply adding a corner square according to some rule?  That would be
an explicit combinatorial (as opposed to algebraic) version of
Sylvester's approach.

Unfortunately we don't know how to construct such a nice
injection.\endnote{The existence of an injection is known and follows
from the \emph{Dilworth theorem} and the \emph{Peck property}, see R.~P.~Stanley,
Quotients of Peck posets, \emph{Order}~\textbf{1} (1984), 29--34. There is
no known efficient algorithm to compute this injection.  Finally, the existence
of a symmetric saturated chain decomposition in this case is a well-known open
problem, see R.~P.~Stanley, Weyl groups, the hard Lefschetz theorem, and the Sperner
property, \emph{SIAM J.~Algebraic Discrete Methods}~\textbf{1} (1980),  168--184. }
It's just the
first of the many frustrations one encounters with algebraic proofs.  Most of them
are simply too rigid to be ``combinatorialized''.   It doesn't mean that there is no
combinatorial interpretation for $C(n,k,\ell)$ at all.  There is one very uninteresting
interpretation due to Panova and myself, based on a very interesting (but cumbersome)
identity by O'Hara.\endnote{A combinatorial interpretation for $C(n,k,\ell)$
based on O'Hara's identity due to I.~Pak and G.~Panova, is given on p.~9 in \ts
\href{http://www.math.ucla.edu/~pak/hidden/papers/Panova_Porto_meeting.pdf}{https://tinyurl.com/ydemhyf5}.
For O'Hara's identity and its simple proof, see I.~G.~Macdonald,
An elementary proof of a $q$-binomial identity, in
{\em $q$-series and partitions}, IMA, Springer, New York, 1989, 73--75.
For O'Hara's combinatorial proof based on the identity, see
D.~Zeilberger, Kathy O'Hara's constructive proof of the unimodality of the
Gaussian polynomials, \emph{Amer.\ Math.\ Monthly}~\textbf{96} (1989), 590–-602. }
Also, from the Computer Science point of view, it is easy
to show that $C(n,k,\ell)$ as a function is in~$\SP$.  We leave it to the reader
to figure out why (or what does that even mean).\endnote{This means that there
is a set of 0-1 sequences whose size is $C(n,k,\ell)$, s.t.\ the membership can
be decided in polynomial time (in~$n$).  This follows immediately from the fact
that $C(n,k,\ell)$ can be computed in polynomial time since Gaussian coefficients
have an explicit recursive formula
$$
\binom{n}{k}_q \, = \, \binom{n-1}{k}_q \. + \, q^{n-k} \binom{n-1}{k-1}_q
$$
and have $O(n^2)$ degree as a polynomial in~$q$. }

To finish this story, we should mention Stanley's 1989 approach to~\eqref{eq:qbin}
using finite group actions.\endnote{R.~P.~Stanley,
Log-concave and unimodal sequences in algebra, combinatorics, and geometry,
in \emph{Ann.\ New York Acad.\ Sci.}~\textbf{576} (1989), 500--535. 
See also White's paper\endnotemark[\ref{en:Kostka}].}\label{en:Stanley}
More recently, Panova and I introduced a different technique based on
properties of the \emph{Kronecker coefficients} of~$S_n$, via the equality
\ts $C(n,k,\ell) = g\bigl((n-k)^k,(n-k)^k,(n-\ell,\ell)\bigr)$.\endnote{I.~Pak, 
G.~Panova, Strict unimodality of $q$-binomial coefficients,
\emph{C.R.\ Acad.\ Sci.\ Paris, Ser.~I.\ Math.}~\textbf{351} (2013), 415--418,
and I.~Pak, G.~Panova, Bounds on certain classes of Kronecker and $q$-binomial
coefficients, \emph{J.\ Combin. Theory~Ser. A.}~\textbf{147} (2017), 1--17.}
Here the Kronecker coefficients $g(\la, \mu,\nu)$ can be defined
as structure constants for products of $S_n$ characters: \.
$\chi^\mu \. \chi^\nu \ts = \ts \sum_\la \. g(\la, \mu,\nu) \. \chi^\la$.\endnote{The
Kronecker coefficients are quite interesting, but very mysterious
and nobody's friend.  They are hard to compute, have no known combinatorial
interpretation, and tend to appear in different areas of science.
For connections to Computational Complexity, see
C.~Ikenmeyer, K.~Mulmuley, M.~Walter,
On vanishing of Kronecker coefficients, \emph{Comput.\ Complexity}~\textbf{26} (2017), 949--992, and
I.~Pak, G.~Panova, On the complexity of computing Kronecker coefficients,
\emph{Comput.\ Complexity}~\textbf{26} (2017), 1--36.
For connections to Quantum Information Theory, see
M.~Christandl, G.~Mitchison,
The spectra of quantum states and the Kronecker coefficients of the symmetric group,
\emph{Comm.\ Math.\ Phys.}~\textbf{261} (2006), 789--797.}
Both approaches imply stronger inequalities than~\eqref{eq:qbin},
but neither gets us closer to a simple injective proof.

\medskip

\nin
We turn now to {\it log-concavity of independent sets:}
\begin{equation}\label{eq:ahk}
a_{k-1}(M)\cdot a_{k+1}(M)\. \le \. a_k(M)^2, \ \ \text{where}
\end{equation}
$a_k(M)$ is the number of independent $k$-subsets of a
matroid~$M$.\endnote{For the background on matroids, see 
N.~White, \emph{Theory of Matroids},  Cambridge Univ.\ Press, 1986, and 
J.~Oxley, \emph{Matroid theory} (Second ed.), Oxford Univ.\ Press, 2011.}
Note that the log-concavity implies unimodality,
and in the special case of a \emph{free matroid} (all elements are independent)
this gives~\eqref{eq:bin}.\endnote{To be clear, for a free matroid $M$ on
$n$ elements, we have $a_k(M)=\binom{n}{k}$.}

The inequality~\eqref{eq:ahk} is a celebrated recent result by Adiprasito,
Huh and Katz (2018), which showed that a certain ``cohomology ring'' associated
with~$M$ satisfies the hard Lefschetz theorem and the Hodge--Riemann relations.
This resolved conjectures by Welsh and Mason (1970s).\endnote{K.~Adiprasito,
J.~Huh, E.~Katz, Hodge theory for combinatorial geometries,
\emph{Ann.\ Math.}~\textbf{188} (2018), 381--452. See also a friendly
introduction by the same authors in 
\emph{Notices AMS}~\textbf{64} (2017), no.~1, 26--30.}

It would be na\"{i}ve for us to ask for a direct combinatorial proof via an injection,
or by some other elementary means.\endnote{This does not mean that no simple
non-combinatorial proof exists.  In recent years, the Adiprasito--Huh--Katz approach
has been substantially simplified and extended.  See
N.~Anari, K.~Liu, S.~Oveis Gharan, C.~Vinzant, Log-Concave Polynomials III:
Mason's Ultra-Log-Concavity Conjecture for Independent Sets of Matroids,
{\tt arXiv:1811.01600}, and P.~Br\"{a}nd\'{e}n, J.~Huh,
Hodge--Riemann relations for Potts model partition functions, {\tt arXiv:1811.01696}. }
For example, Stanley in 1981 used the
Aleksandrov--Fenchel inequalities in convex geometry to prove that the
log-concavity is preserved under taking truncated sum with a
free matroid\endnote{To understand the abstract definition of the
truncated sum of matroids,
think of two sets of vectors $S$ and~$S'$ in vector spaces $V$ and~$V'$,
of dimensions $d$ and $d'$, respectively. Suppose $\max\{d,d'\}\le k\le d+d'$.
The truncated matroid is given by projection $S\times S'\ssu V\times V'$
onto a generic $k$-subspace.},
already an interesting but difficult special case proved by inherently
non-combinatorial means.\endnote{R.~P.~Stanley, Two combinatorial
applications of the Aleksandrov--Fenchel inequalities, \emph{J.~Combin.
Theory~Ser.~A}~\textbf{31} (1981), 56--65.
For the extension of this proof to truncated sums of general matroids
(using a different terminology and motivation), see L.~Gurvits, A short proof,
based on mixed volumes, of Liggett's theorem on the convolution of
ultra-logconcave sequences, \emph{El.~J.\ Combin.}~\textbf{16} (2009),
no.~1, Note~5, 5~pp.}\label{en:Gurv}

There is also a Computational Complexity version of the
problem which might be of interest.  Let \ts
$A(k,M):= a_k(M)^2 - a_{k-1}(M)\cdot a_{k+1}(M)$.  Does $A(k,M)$ count any set of
combinatorial objects?

For the sake of clarity, let $G=(V,E)$ be a simple connected graph and $M$
the corresponding matroid, i.e.\ bases in $M$ are spanning trees in $G$.
Then $a_k(M)$ is the number of spanning forests in~$G$ with $k$~edges.
Note that computing $a_k(M)$ is $\SP$-complete in full
generality.\endnote{For a comprehensive discussion of complexity of this
and other values of the Tutte polynomial, see
D.~J.~A.~Welsh, \emph{Complexity: knots, colourings and counting},
Cambridge Univ.\ Press, 1993.}
Therefore, computing $A(k,M)$ is $\SP$-hard.

Now, $A(k,M)$ is in $\GapP$, i.e.\ equal to the difference of two $\SP$-functions.
Does $A(k,M)$ lie in~$\SP$?  This seems unlikely, but the current state
of art of Computational Complexity doesn't seem to provide us with tools
to even approach a negative solution.\endnote{For the recent overview of
complexity of computing combinatorial sequences, see I.~Pak, Complexity
problems in enumerative combinatorics, in \emph{Proc.\ ICM Rio de Janeiro}, Vol.~3,
2018, 3139--3166; revised and expanded version in {\tt arXiv:1803.06636}.}

\medskip

\nin
To fully appreciate the last example, consider the
{\it log-concavity of matching numbers:}
\begin{equation}\label{eq:match}
m_{k-1}(G)\cdot m_{k+1}(G)\. \le \. m_k(G)^2, \ \ \text{where}
\end{equation}
$m_k(G)$ is the number of $k$-matchings in a simple graph $G=(V,E)$,
i.e.\ $k$-subsets of edges which are pairwise disjoint.  For example, \ts
$m_n(K_{2n})=(2n-1)\cdots 3 \cdot 1$.   While perfect matchings don't necessarily
define a matroid, they do have a similar flavor from a Combinatorial
Optimization point of view.\endnote{By that we mean that
the linear programming works in polynomial time on both perfect
matchings in graphs and bases of matroids. For an overview and
a careful explanation of the connection, see W.~H.~Cunningham,
Matching, matroids, and extensions,
\emph{Math.\ Program.}~\textbf{91} (2002), Ser.~B, 515--542.}
The inequality~\eqref{eq:match} goes back
to Heilmann and Lieb~(1972) and is a rare case when the injection
strategy works well.\endnote{O.~J.~Heilmann, E.~H.~Lieb,
Theory of monomer-dimer systems, \emph{Comm.\ Math.\ Phys.}~\textbf{25}
(1972), 190--243.}  The following argument is due to
Krattenthaler (1996).\endnote{C.~Krattenthaler,
Combinatorial proof of the log-concavity of the sequence
of matching numbers,
\emph{J.~Combin.\ Theory Ser.~A}~\textbf{74} (1996), 351--354.}

Take a $(k-1)$-matching $\be$ whose edges we color \emph{blue} and a
$(k+1)$-matching $\ga$ whose edges we color \emph{green}.  The union
$\be\cup \ga$ of these two sets of edges splits into connected
components, which are either paths or cycles, all alternately
colored.  Ignore for the time being all cycles and paths of even lengths.
Denote by $(r-1)$ the number of odd-length paths which have extra color blue.
There are then $(r+1)$ odd-length paths which have extra color green.

Now, allow switching colors in any of the $2r$ odd-length paths.  After
recoloring, we want to have $r$ odd-length paths extra color blue and
the same with green.  This amounts to a constructive injection
from $(r-1)$-subsets of $[2r]$ to $r$-subsets of $[2r]$, which we already
know how to do as a special case of proving~\eqref{eq:bin}.

We leave to the reader the problem of finding an explicit combinatorial interpretation
for \. $M(k,G):=m_k(G)^2-m_{k-1}(G)\cdot m_{k+1}(G)$, proving that this function is
in~$\SP$.  Note that computing $m_k(G)$ is famously $\SP$-complete, which
implies that so is $M(k,G)$.  This
makes the whole connection to Computational Complexity even more confusing.
What exactly makes matchings special enough for this argument
to work?\endnote{For several explicit combinatorial injections in the context of
log-concavity, see B.~E.~Sagan, Inductive and injective proofs of log concavity results,
\emph{Discrete Math.}~\textbf{68} (1988), 281--292.  In the opposite direction,
see an example of a log-covex combinatorial sequence in S.~DeSalvo, I.~Pak,
Log-concavity of the partition function, \emph{Ramanujan J.}~\textbf{38}
(2015), 61--73. Here the result is that partition
function satisfies \ts $p(n)^2-p(n-1)\ts p(n+1)\ge 0$, for all $n>25$.
It is unlikely this inequality has a direct combinatorial proof. }

\medskip

\nin
If there is any pattern to the previous examples, it can be summarized as follows:
the deeper one goes in an algebraic direction, the more involved are the
inequalities and the less of a chance of a combinatorial proof.  To underscore
this point, consider the following three \emph{Young tableaux
inequalities}:
\begin{equation}\label{eq:yt}
\bigl(f^\la\bigr)^2 \. \le \. n!, \ \quad \bigl(c^\la_{\mu\ts\nu}\bigr)^2 \.
\le \. \binom{n}{k}, \ \quad c^\la_{\mu\ts \nu} \. \le \. c^\la_{\mu\vee\nu,\ts \mu \wedge \nu}\,, \
\quad \text{for all} \ \
\la\vdash n, \. \mu \vdash k, \. \nu \vdash n-k.
\end{equation}
Here \. $f^\la = \bigl|\text{SYT}(\la)\bigr|$ \. is the number of standard Young tableaux
of shape~$\la$, equal to the dimension of the corresponding irreducible $S_n$-module as above.
Similarly, \ts $c^\la_{\mu\ts \nu}= \bigl|\text{LR}(\la/\mu,\nu)\bigr|$ \. is the
\emph{Littlewood--Richardson coefficient}, equal to the number of Littlewood--Richardson
tableaux of shape $\la/\mu$ and weight~$\nu$.  It can be defined as a structure constant
for products of Schur functions: \. $s_\mu \. s_\nu \. = \. \sum_\la \. c^\la_{\mu\ts \nu} \. s_\la$.
Finally, $\mu\vee\nu$ and $\mu \wedge \nu$ denote the union and intersection, respectively,
of the corresponding Young diagrams.\endnote{For a comprehensive introduction to
combinatorics of Young tableaux and the Littlewood--Richardson coefficients, see
B.~E. Sagan, {\em The Symmetric Group}, Springer, New York, 2001, and
R.~P.~Stanley, \emph{Enumerative Combinatorics}, Vol.~2,
Cambridge U.~Press, Cambridge, 1999, Ch.~7.}

Now, the first inequality in~\eqref{eq:yt} is trivial algebraically, but
its combinatorial proof is highly nontrivial -- it is a restriction of
the \emph{RSK correspondence}.\endnote{The RSK is a fundamental bijection
described in both Sagan and Stanley's books, see above.}
The second inequality is quite recent and follows easily from the definition
and the Frobenius reciprocity.  We believe it is unlikely that there is
a combinatorial injection, even though there is a nice double
counting argument.\endnote{The inequality is given in
I.~Pak, G.~Panova, D.~Yeliussizov,
On the largest Kronecker and Littlewood--Richardson coefficients,
\emph{J.~Combin.\ Theory Ser.~A}~\textbf{165} (2019), 44--77, $\S$4.
The ``double counting injection'' is given by the following one line argument:
$$\sum_{\la \vdash n} \.  \bigl(c^\la_{\mu,\ts \nu}\bigr)^2 \ts f^\mu\ts f^\nu
\, = \, \sum_{\la \vdash n} \.  c^\la_{\mu,\ts \nu} \ts \bigl(c^\la_{\mu,\ts \nu} \ts f^\mu\ts f^\nu\bigr)
\, \le \, \sum_{\la \vdash n} \.
c^\la_{\mu,\ts \nu} \. f^\la \, = \,  \binom{n}{k} \ts f^\mu \ts f^\nu.
$$
Here the inequality and the last equality can be made injective by
using two different combinatorial interpretations of the
Littlewood--Richardson coefficients.
}

Finally, the third inequality in~\eqref{eq:yt} is a corollary of the powerful inequality
by Lam, Postnikov and Pylyavskyy (2007) using the curious Temperley--Lieb immanant
machinery.\endnote{T.~Lam, A.~Postnikov, P.~Pylyavskyy, Schur positivity and
Schur log-concavity, \emph{Amer.~J.\ Math.}~\textbf{129} (2007), 1611--1622.
The inequality there is more general and stated in a different but equivalent way;
to see the connection, consider only straight shapes and write the inequalities
in Schur functions basis. Also, this inequality is a refinement of an 
older, FKG-type inequality: 
$$
f^{\la}\. \cdot \. f^{\mu} \. \le \. f^{\la\vee \mu}\. \cdot \. f^{\la\wedge \mu}. 
$$
The proof of the latter is elementary and based on the hook-length formula, 
see A.~Bj\"{o}rner, A $q$-analogue of the FKG inequality and some applications, 
\emph{Combinatorica}~\textbf{31} (2011), 151--164. }
The key ingredient in the proof is Haiman's theorem which in turn uses the
Kazhdan--Lusztig conjecture proven by Beilinson--Bernstein and
Brylinski--Kashiwara.\endnote{M.~Haiman, Hecke algebra characters and immanant
conjectures, \emph{JAMS}~\textbf{6} (1993), 569--595.} While stranger things
have happened, we would be very surprised if this inequality had a simple
combinatorial proof.\endnote{Another notable example is the inequality
implied by the \emph{Naruse hook-length formula} (NHLF):
$$(\diamondsuit) \qquad f^{\la/\mu} \, \ge \, n! \. \prod_{(i,j)\in \la/\mu} \. \frac{1}{h(i,j)}\,,
$$
where \ts $f^{\la/\mu} = \bigl|\SYT(\la/\mu)\bigr|$ \ts is the number of standard
Young tableaux of skew shape and $h(i,j) = \la_i+\la_j'-i-j+1$ is the hook-length in~$\la$.
When $\mu=\emp$, the inequality becomes an equality; it is the celebrated
\emph{hook-length formula} with a plethora of bijective, probabilistic, and
algebraic proofs.  Unfortunately, no direct injective proof is known for~$(\diamondsuit)$.
In a special case of $\la/\mu=\tau^\ast$, a Young diagram~$\tau$ rotated 180$^\circ$,
it gives the following \emph{hook inequality}:
$$(\heartsuit) \qquad \prod_{(i,j)\in \tau} \. h(i,j) \, \le \, \prod_{(i,j)\in \tau} \. h^\ast(i,j),
$$
where $h^\ast(i,j)=i+j-1$ are the hooks in $\tau^\ast$.  Even this inequality
does not seem to have a direct injective proof.  For $(\heartsuit)$ and
further discussion and references to the NHLF, see A.~Morales, I.~Pak, G.~Panova,
Asymptotics of the number of standard Young tableaux of skew shape,
\emph{Eur.\ J.\ Combin.}~\textbf{70} (2018), 26--49.  For a recent combinatorial
(but not injective!) proof of $(\heartsuit)$, see I.~Pak, F.~Petrov, V.~Sokolov,
Hook inequalities, \ts {\tt arXiv:1903.11828}. }

\medskip

\nin
We conclude on a positive note, with a combinatorial inequality
where everything works as well as it possibly could.  Consider the
following \emph{majorization property of contingency tables}:
\begin{equation}\label{eq:ct}
\rT(\bba, \bbb) \, \le \,  \rT(\bba', \bbb') \ \ \. \text{for all} \ \
\bba' \trianglelefteq \bba, \, \bbb' \trianglelefteq \bbb.
\end{equation}
Here \. $\mathbf a=(a_1,\dots,a_m)$, $a_1\ge \ldots \ge a_m>0$,
and $\mathbf b=(b_1,\dots,b_n)$, $b_1\ge \ldots \ge b_n>0$, are
two integer sequences with equal sum:
$$
\sum_{i=1}^m \ts a_i \, = \, \sum_{j=1}^n \ts b_j \. = \. N.
$$
A \emph{contingency table} with margins $(\mathbf a, \mathbf b)$
is an $m \times n$ matrix of non-negative integers whose $i$-th row sums to $a_i$
and whose $j$-th column sums to~$b_j$, for all $i\in [m]$ and $j\in [n]$.
$\rT(\bba, \bbb)$ denotes the number of all such matrices.
Finally, for sequences $\bba$ and $\bba'$ with the same sum, we write \ts
$\bba \trianglelefteq \bba'$ \ts if \ts $a_1\le a_1'$, $a_1+a_2\le a_1'+a_2'$,
$a_1+a_2+a_3\le a_1'+a_2'+a_3'$, $\ldots$ \ts  In other words, the
inequality~\eqref{eq:ct} says that there are more contingency tables
when the margins are more evenly distributed.

Contingency tables can be viewed as adjacency matrices of
bipartite multi-graphs with given degree
distribution.  They play an important role in Statistics
and Network Theory.\endnote{See e.g.\ P.~Diaconis, A.~Gangolli,
		Rectangular arrays with fixed margins,
		\emph{Disc.\ Prob.\ Alg.}~\textbf{72} (1995), 15--41, and
B.~S.~Everitt, \emph{The analysis of contingency tables} (Second ed.),
Chapman \& Hall, London, 1992.}
We learned the inequality~\eqref{eq:ct} from a paper by Barvinok (2007),
but it feels like something that should have been known for
decades.\endnote{See Eq.~(4) on p.~111 in  A.~Barvinok,
Brunn--Minkowski inequalities for contingency tables and integer flows,
\emph{Adv.\ Math.}~\textbf{211} (2007), 105--122. Note that a
more general inequality~(3) in that paper remains open. }

Now, we know two fundamentally different proofs of~\eqref{eq:ct}.  The first
is an algebraic proof using Schur functions which amounts to proving
the following standard inequality for Kostka numbers: \.
$K_{\la\ts\mu}\le K_{\la\ts\nu}$ \. for all \. $\mu \trianglerighteq \nu$,
where $K_{\la\ts\mu}$ is the number of semistandard Young tableaux of shape~$\la$
and weight~$\mu$.
This inequality can also be proved directly, so combined with the RSK we obtain
an injective proof of~\eqref{eq:ct}.\endnote{Barvinok's proof is exactly the
combinatorial identity given by the RSK written in terms of the Kostka numbers.
For a direct injective proof of the inequality for Kostka numbers, see
D.~E.~White, Monotonicity and unimodality of the pattern inventory, 
\emph{Adv.\ Math.}~\textbf{38} (1980), 101--108.}\label{en:Kostka}

Alternatively, one can prove the inequality directly for $2\times n$
rectangles and $(+1,-1)$ changes in row (column) sums.\endnote{To see this,
note that for $\bba=(a_1,\ldots,a_n)$, $\bbb = (N-k,k)$, $N=a_1+\ldots+a_n$,
we have:
$$
\rT(\bba, \bbb) \, = \, \bigl[q^k\bigr] \, \prod_{i=1}^n \. \bigl(1+q+\ldots + q^{a_i-1}\bigr)
$$
Each term in the product is unimodal, thus so is the product
(see e.g.\ Stanley's paper\endnotemark[\ref{en:Stanley}] above). This
implies unimodality.  To get a direct injection, see the Greene--Kleitman
symmetric chain decomposition approach described above.\endnotemark[\ref{en:GK}] 
Curiously, White's argument\endnotemark[\ref{en:Kostka}] 
can be also translated into this language. }\label{en:two-rows}
%
Combining these injections
together gives a cumbersome, yet explicit injection.  In principle, either
of the two approaches can then be used to give a combinatorial interpretation
for \.  $\rT(\bba', \bbb')-\rT(\bba, \bbb)$.\endnote{I am being very vague
here both in an effort to avoid technicalities and at least to some extend
not to reveal the work in progress.  If and when this work is completed, 
I will update this note.}

\smallskip

In conclusion, let us note that we came full circle.  Let \ts
$m=2$, $a_1=n-k+1$, $a_2=k-1$, $a_1'=n-k$, $a_2'=k$, and \ts
$b_1=\ldots=b_n=b_1'=\ldots=b_n'=1$.\endnote{This is the case in the
 note above\endnotemark[\ref{en:two-rows}], of course.} Observe that \ts
$\rT(\bba, \bbb)=\binom{n}{k-1}$ \ts and \ts
$\rT(\bba', \bbb')=\binom{n}{k}$.  The inequality~\eqref{eq:bin} is a
special case of~\eqref{eq:ct} then.

\vskip.5cm

\subsection*{Acknowledgements}
We are grateful to Karim Adiprasito, Sasha Barvinok, Sam Dittmer,
Tom Liggett, Alex Mennen, Alejandro Morales, Greta Panova, F\"edor Petrov,
Richard Stanley and Damir Yeliussizov for many interesting and
helpful conversations.  The author was partially supported by the NSF.


\newpage

\printendnotes*

\end{document}